\documentclass{amsart}
\usepackage{graphicx}
\newtheorem{theorem}{Theorem}[section]
\newtheorem{lemma}[theorem]{Lemma}

\newtheorem{proposition}[theorem]{Proposition}
\usepackage{graphicx}

\theoremstyle{definition}

\theoremstyle{remark}

\numberwithin{equation}{section}

\begin{document}

\markboth{Masakazu Teragaito}
{Hyperbolic knots with three toroidal Dehn surgeries}


\title{Hyperbolic knots with three toroidal Dehn surgeries}

\author{Masakazu Teragaito}
\address{Department of Mathematics and Mathematics Education,
Hiroshima University, 1-1-1 Kagamiyama, Higashi-hiroshima 739-8524, Japan}

\begin{abstract}
It is conjectured that a hyperbolic knot admits at most three Dehn surgeries
which yield closed $3$-manifolds containing incompressible tori.
We show that there exist infinitely many hyperbolic knots
which attain the conjectural maximum number.
Interestingly, those surgeries correspond to consecutive integers.
\end{abstract}



\maketitle

\section{Introduction}\label{sec:intro}

For a knot $K$ in the $3$-sphere $S^3$, let $E(K)=S^3-\mathrm{Int}N(K)$ be its exterior, where $N(K)$ denotes a tubular neighborhood of $K$.
A \textit{slope\/} $\alpha$ is the isotopy class of an unoriented essential simple loop on $\partial E(K)$.
Slopes are parameterized by the set $\mathbb{Q}\cup \{1/0\}$ in the usual way (see \cite{Ro}).
In particular, $1/0$ corresponds to the slope of a meridian.
A slope is said to be \textit{integral\/} if it corresponds to an integer.
Thus an integral slope runs once along the knot.
For a slope $\alpha$, let $K(\alpha)$ be the closed orientable $3$-manifold obtained by \textit{$\alpha$-Dehn surgery\/} on $K$,
that is, $K(\alpha)$ is the union of $E(K)$ and a solid torus $V$, where $V$ is attached to $E(K)$ along their boundaries so that
a meridian of $V$ goes to $\alpha$ on $\partial E(K)$.

Assume that $K$ is a hyperbolic knot.
When $K(\alpha)$ is not hyperbolic, the slope $\alpha$ is called an \textit{exceptional slope}, 
and the surgery is also said to be \textit{exceptional}.
Each hyperbolic knot has only finitely many exceptional slopes by Thurston's hyperbolic Dehn surgery theorem \cite{Th}.
It is conjectured that except the figure-eight knot and the $(-2,3,7)$-pretzel knot, any hyperbolic knot
admits at most six exceptional slopes \cite[Problem 1.77(A)(1)]{Ki}.
On the other hand, the resulting manifold by an exceptional Dehn surgery is expected to be
either $S^3$, a lens space, a Seifert fibered manifold over the $2$-sphere with three exceptional fibers
 (referred to as a small Seifert fibered manifold), or a toroidal manifold \cite{Go}.
Here, a \textit{toroidal manifold\/} is a closed $3$-manifold which contains an incompressible torus.

According to the type of the resulting manifold, an exceptional surgery is referred to as 
a \textit{lens space surgery\/}, a \textit{Seifert surgery\/} or a \textit{toroidal surgery\/}, respectively.
By Gordon and Luecke's theorem \cite{GL}, only $1/0$-Dehn surgery can yield $S^3$ for a non-trivial knot.
The cyclic surgery theorem \cite{CGLS} implies that a hyperbolic knot admits at most two lens space surgeries, which must be integral, and 
if there are two, they are consecutive.
In fact, there are infinitely many hyperbolic knots with two lens space surgeries.
Except the figure-eight knot with six Seifert surgeries, a hyperbolic knot seems to admit at most three Seifert surgeries.
Recently, Deruelle, Miyazaki and Motegi \cite{DMM} gave a hyperbolic knot with three Seifert surgeries corresponding to
any successive three integers.

In this paper, we will focus on toroidal surgeries.
Eudave-Mu\~{n}oz \cite{EM} conjectured that any hyperbolic knot admits at most three toroidal surgeries (see also \cite[Problem 1.77(A)(5)]{Ki}).
This conjecture holds for $2$-bridge knots \cite{BW} and Montesinons knot \cite{W}.
In general, the best result forward this direction is Gordon and Wu's one \cite{GW} which claims that
a hyperbolic knot admits at most four toroidal surgeries, and if there are four, then they correspond to consecutive integers. 
As far as we know, the only examples of hyperbolic knots that realize the expected optimum number are
the figure-eight knot and the $(-2,3,7)$-pretzel knot, with toroidal slopes $\{-4,0,4\}$ and $\{16,37/2,20\}$, respectively.
The purpose of this paper is to give the first infinite family of hyperbolic knots with three toroidal surgeries.
Interestingly, these toroidal surgeries correspond to consecutive integers. 

\begin{theorem}\label{thm:main}
There are infinitely many tunnel number one, hyperbolic knots, each of which
admits three toroidal Dehn surgeries corresponding to consecutive integers.
\end{theorem}

Our construction is based on the Montesinos trick \cite{M}.
We will construct a tangle which can produce the unknot by summing some rational tangle.
This implies that the double branched cover of the tangle gives the exterior of a knot in $S^3$.
The tangle is carefully given so that the sums with three rational tangles yield knots or links which admit
essential Conway spheres.  Thus our knot will admit three toroidal surgeries.
The idea of the tangle is a variation of the pentangle introduced by Gordon and Luecke in \cite{GL1} for a different purpose.
Although it is easy to see that our tangle admits two rational tangle sums yielding essential Conway spheres,
it came as a surprise that the third rational tangle sum also yields an essential Conway sphere.

As an additional interesting feature, one of the toroidal surgeries for each of our knots yields a closed $3$-manifold which contains
a unique incompressible torus meeting the core of the attached solid torus in exactly four points, but does not contain
an incompressible torus meeting the core in less than four points.
Such examples of toroidal surgeries on hyperbolic knots have been already given by
Eudave-Mu\~{n}oz \cite{EM2}.
The simplest knot among his knots seems to have genus $37$ as he wrote, but
our simplest knot, as shown in Fig.~\ref{fig:k2}, has genus nine.

\section{Construction}\label{sec:const}

A \textit{tangle\/} is a pair $(B,t)$ where $B$ is a $3$-ball and $t$ is a finite disjoint union
of simple closed curves and properly embedded arcs.

\begin{figure}[th]
\centerline{
\includegraphics*[scale=0.4]{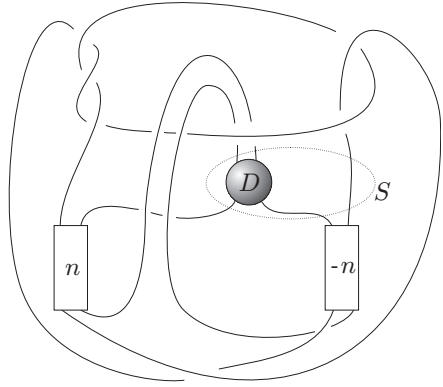}}
\vspace*{8pt}
\caption{The tangle $T_n$}\label{fig:tangle}
\end{figure}

For an integer $n\ge 2$, consider the tangle $T_n=(B,t)$
as shown in Fig.~\ref{fig:tangle}, where $B$ is the $3$-ball obtained from $S^3$ by removing
the interior of the $3$-ball $D$.
Here, the rectangle labeled by an integer $n$ (resp. $-n$)
denotes $n$ right-handed (resp. left-handed) vertical half-twists.
(Although $T_n$ can be defined for any integer $n$, the restriction $n\ge 2$ suffices to prove our result.)

We will insert several rational tangles into $D$, giving a knot or link in $S^3$.
In fact, we use only the rational tangles illustrated in Fig.~\ref{fig:tangle}, where
we adopt the convention of \cite{EM3}.

\begin{figure}[th]
\centerline{
\includegraphics*[scale=0.4]{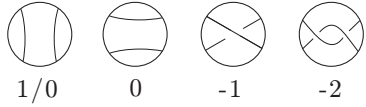}}
\vspace*{8pt}
\caption{Rational tangles}\label{fig:ratio}
\end{figure}

A filling of $T_n$, $T_n(\alpha)$, refers to filling $D$ with the rational tangle of slope $\alpha$.
Let $\widetilde{T}_n(\alpha)$ denote the double branched cover of $S^3$ branched over $T_n(\alpha)$.

\begin{lemma} $\widetilde{T}_n(1/0)=S^3$.
\end{lemma}

\begin{proof}
This easily follows from the figures.
After inserting the $1/0$-tangle into $D$, the two twist boxes are canceled, and
we see that $T_n(1/0)$ is the unknot.
Thus the double branched cover  $\widetilde{T}_n(1/0)$ is $S^3$.
\end{proof}

Since the lift of a rational tangle is a solid torus,
the lift of the $3$-ball $B$ of the tangle $T_n=(B,t)$
gives the exterior of a knot in $S^3$, which is denoted by $K_n$.

We use $M(r,s)$ to denote the Montesinos tangle consisting of two rational tangles
corresponding to the rational numbers $r$ and $s$, respectively.
(See \cite{GL2}.)
The double branched cover of $M(r,s)$, denoted by $D^2(r,s)$, is the Seifert fibered manifold
over the disk with two exceptional fibers of type $r$ and $s$. 

\begin{lemma}\label{lem:slope0}
$\widetilde{T}_n(0)=D^2(1/2,1/3)\cup D^2(1/n,-1/(n+1))$.
\end{lemma}

\begin{figure}[th]
\centerline{
\includegraphics*[scale=0.7]{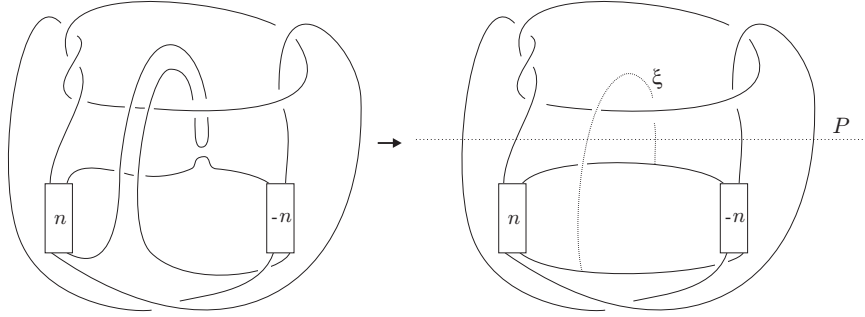}}
\vspace*{8pt}
\caption{$0$-filling}\label{fig:slope0}
\end{figure}

\begin{proof}
Figure \ref{fig:slope0} shows $T_n(0)$ which is decomposed along a tangle sphere $P$ into
two Montesions tangles $M(1/2,1/3)$ and $M(1/n,-1/(n+1))$.
Thus $\widetilde{T}_n(0)$ is decomposed along a torus into two Seifert fibered manifolds
$D^2(1/2,1/3)$ and $D^2(1/n,-1/(n+1))$, where the Seifert fibers intersect once on the torus.
\end{proof}

\begin{lemma}\label{lem:slope-1}
$\widetilde{T}_n(-1)=D^2(1/2,1/n)\cup D^2(-1/2,-1/(n+1))$.
\end{lemma}

\begin{figure}[th]
\centerline{
\includegraphics*[scale=0.7]{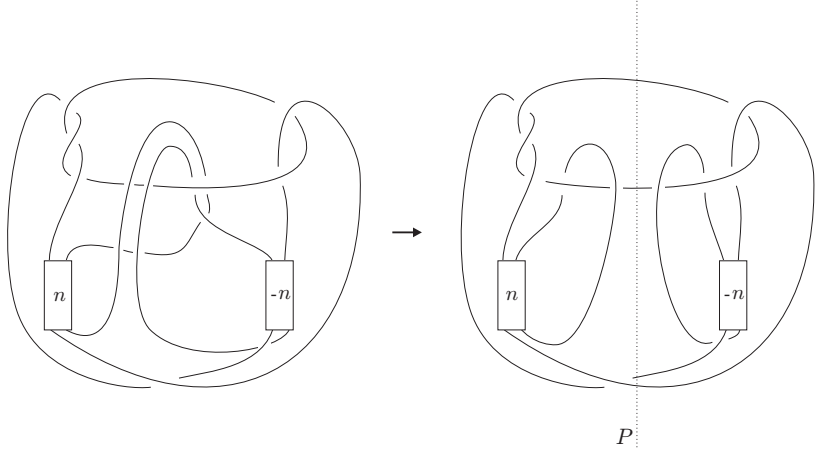}}
\vspace*{8pt}
\caption{$(-1)$-filling}\label{fig:slope-1}
\end{figure}

\begin{proof}
It follows from Fig.~\ref{fig:slope-1} similar to the proof of Lemma \ref{lem:slope0}. 
\end{proof}

\begin{lemma}\label{lem:slope-2}
$\widetilde{T}_n(-2)=D^2(-2/3,1/(n+1))\cup D^2(-2/3,-1/n)$.
\end{lemma}

\begin{figure}[th]
\centerline{
\includegraphics*[scale=0.68]{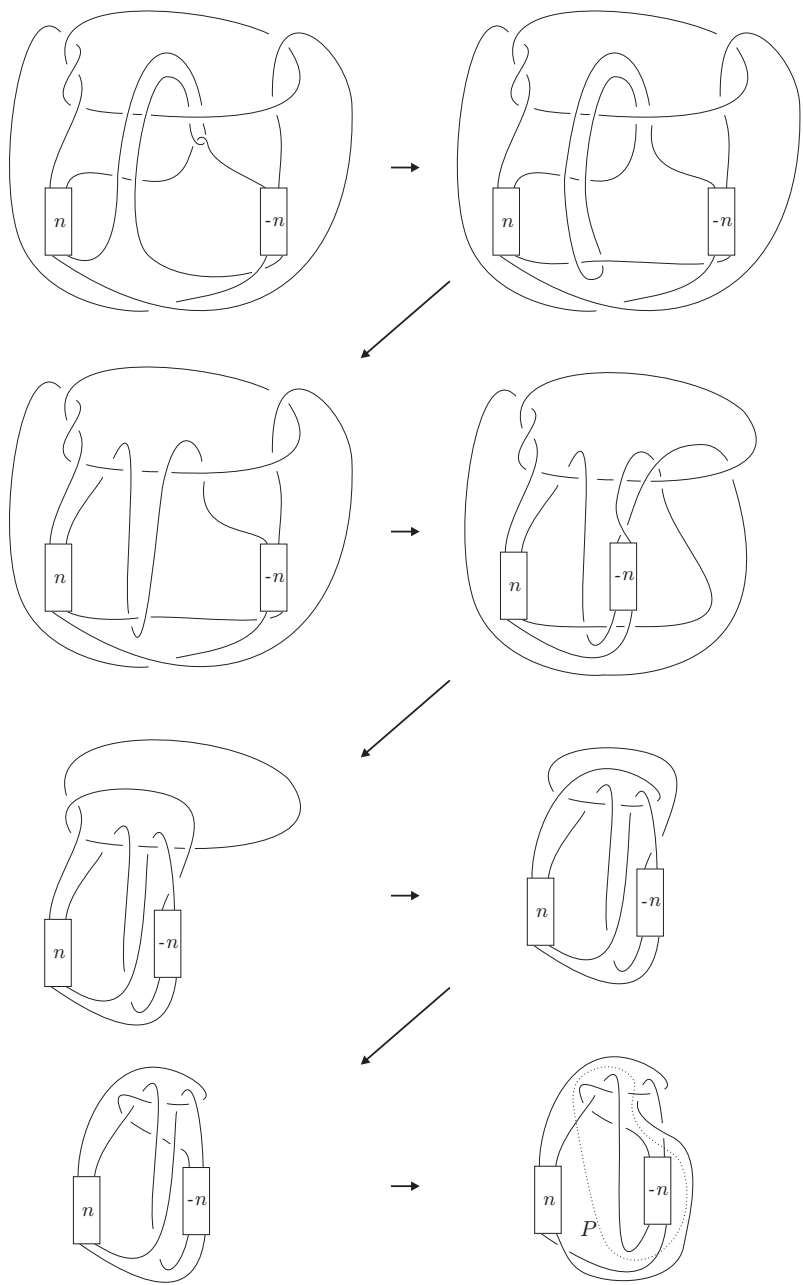}}
\vspace*{8pt}
\caption{$(-2)$-filling}\label{fig:slope-2}
\end{figure}

\begin{proof}
Following the sequence of isotopies as in Fig.~\ref{fig:slope-2},
$T_n(-2)$ is decomposed along a tangle sphere $P$ into two Montesinos tangles as desired.
\end{proof}

\section{Properties of $K_n$}\label{sec:property}

\begin{lemma}\label{lem:3tor}
$K_n$ admits three toroidal slopes which correspond to three successive integers.
Moreover, one of the toroidal surgeries yield a closed $3$-manifold which contains
an incompressible torus meeting the core of the attached solid torus in four points, but does not contain
an incompressible torus meeting the core in less than four points.
\end{lemma}

Those integral slopes will be calculated in the next section by using an explicit description of $K_n$.

\begin{proof}
The first conclusion immediately follows from Lemmas \ref{lem:slope0}, \ref{lem:slope-1} and \ref{lem:slope-2}.
We remark that if $0$-filling for $T_n$ lifts to an integral slope $m$ for $K_n$, then
$(-1)$-filling and $(-2)$-filling lift to $m+1$ and $m+2$, respectively.

In the double branched cover $\widetilde{T}_n(0)$, the lift of the arc $\xi$ as shown in Fig.~\ref{fig:slope0} gives
the core $c$ of the attached solid torus of the surgery on $K_n$ corresponding to the $0$-filling $T_n(0)$.
Since $\xi$ meets the Conway sphere $P$ in two points,
$c$ meets the incompressible torus $R$ as the lift of $P$ in four points.
Let $M_1$ and $M_2$ be the Montesions tangles in the tangle decomposition of $\widetilde{T}_n(0)$ as shown in Fig.~\ref{fig:slope0},
and let $F_i$ be a disk in $M_i$ which divides $M_i$ into two rational tangles for $i=1,2$.
Moreover, we can choose $F_i$ so that $\xi\cap M_i\subset F_i$ for $i=1,2$.
Then each component of $\xi\cap M_i$ is either an arc going from $\partial F_i$ to an intersection point with the strings, or
a spanning arc in $F_i$ which splits it into two disks, each having a point of intersection with the strings.
This implies that $c$ and $R$ intersect minimally in four points by \cite[Example 1.4]{EM2}.

It is well known that each $M_i$ admits the unique Seifert fibration \cite{Wa}.
Since  the Seifert fibers of each side intersect once on the torus $R$ by Lemma \ref{lem:slope0},
$\widetilde{T}_n(0)$ does not admit a Seifert fibration.
In other words, $\{R\}$ gives the torus decomposition of $\widetilde{T}_n(0)$
in the sense of Jaco-Shalen \cite{JS} and Johannson \cite{Jo}.
Thus $\widetilde{T}_n(0)$ contains a unique incompressible torus, which implies the second conclusion.
\end{proof}

We remark that for the other two toroidal surgeries for $K_n$, the resulting manifold contains an incompressible torus
which meets the core of the attached solid torus in two points.

Recall that a knot $K$ has \textit{tunnel number one\/} if there exists an arc $\tau$ with $K\cap \tau=\partial \tau$ such that
$S^3-\mathrm{Int}N(K\cup\tau)$ is a genus two handlebody.
Then such an arc $\tau$ is called an \textit{unknotting tunnel\/} for $K$.

\begin{lemma}\label{lem:tunnel}
$K_n$ has tunnel number one.
\end{lemma}

\begin{figure}[th]
\centerline{
\includegraphics*[scale=0.5]{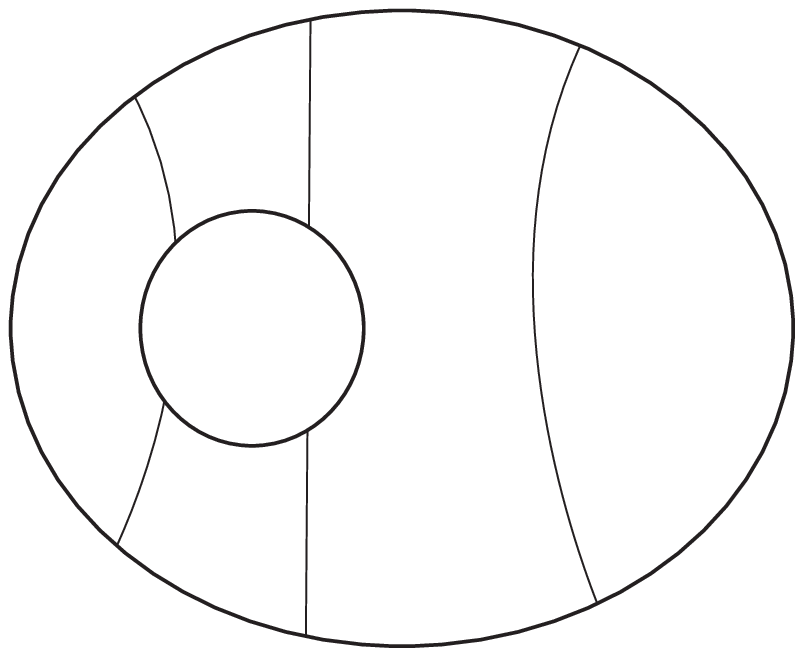}}
\vspace*{8pt}
\caption{}\label{fig:3bridge}
\end{figure}

\begin{proof}
The sphere $S$ illustrated in Fig.~\ref{fig:tangle} splits the tangle $T_n$ into two parts,
one being a $3$-string trivial tangle, and the other being as shown in Fig.~\ref{fig:3bridge}.
(Such a decomposition is called a $3$-bridge decomposition in \cite{EM2}.)
The lift of this decomposition to the exterior of $K_n$ gives a genus two Heegaard splitting.
Thus $K_n$ has tunnel number at most one.
Since $K_n$ admits a toroidal surgery by Lemma \ref{lem:3tor}, $K_n$ is non-trivial, so its tunnel number is one.
\end{proof}

In Fig.~\ref{fig:final}, an unknotting tunnel $\tau$ for $K_n$ is shown by a dotted line.

\begin{lemma}\label{lem:hyp}
$K_n$ is hyperbolic.
\end{lemma}

\begin{proof}
The argument is the same as the proof of Theorem 3.2(3) in \cite{EM2}.
Assume that $K_n$ is not hyperbolic.  Then $K_n$ is either a torus knot or a satellite knot.
Since no surgery on a torus knot produces an incompressible separating torus, $K_n$ is not a torus knot.
Hence $K_n$ is a satellite knot.
Because $K_n$ has tunnel number one by Lemma \ref{lem:tunnel},
$K_n$ is a satellite of some torus knot by \cite{EM1, MS}.
Let $Q$ be the incompressible torus in $E(K_n)$ which bounds the torus knot exterior.
More precisely, $Q$ decomposes $S^3$ into $V\cup N$, where $N$ is the torus knot exterior and $V$ is a solid torus containing
$K_n$ in its interior.
Let $\alpha$ be the toroidal slope corresponding to $0$-filling of $T_n$.
Since the toroidal manifold $K_n(\alpha)=\widetilde{T}_n(0)$ contains the unique incompressible separating torus which
meets the core of the attached solid torus in four points by Lemmas \ref{lem:slope0} and \ref{lem:3tor}, $Q$ is compressible in $K_n(\alpha)$.
This means that the boundary torus $\partial V$ of $V$ is compressible after performing $\alpha$-surgery along $K_n$.
By \cite{S}, the resulting manifold obtained from $V$ by $\alpha$-surgery on $K_n$ is
either a solid torus or the connected sum of a solid torus and a lens space.
The latter is impossible, because $K_n(\alpha)$ is irreducible and not a lens space.
If the former happens, then $K_n(\alpha)$ is obtained by surgery along a torus knot.
This contradicts the fact that any surgery on a torus knot does not produce an incompressible separating torus \cite{Mo}.
\end{proof}

\section{Explicit description of $K_n$}\label{sec:knot}

\begin{figure}[th]
\centerline{
\includegraphics*[scale=0.7]{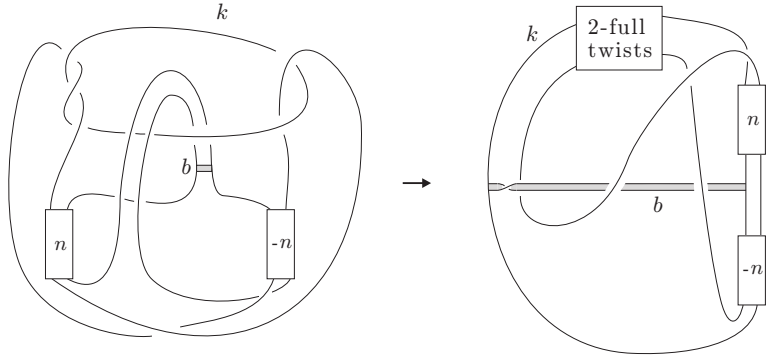}}
\vspace*{8pt}
\caption{The unknot $k$ with band $b$}\label{fig:start}
\end{figure}

\begin{figure}[th]
\centerline{
\includegraphics*[scale=0.7]{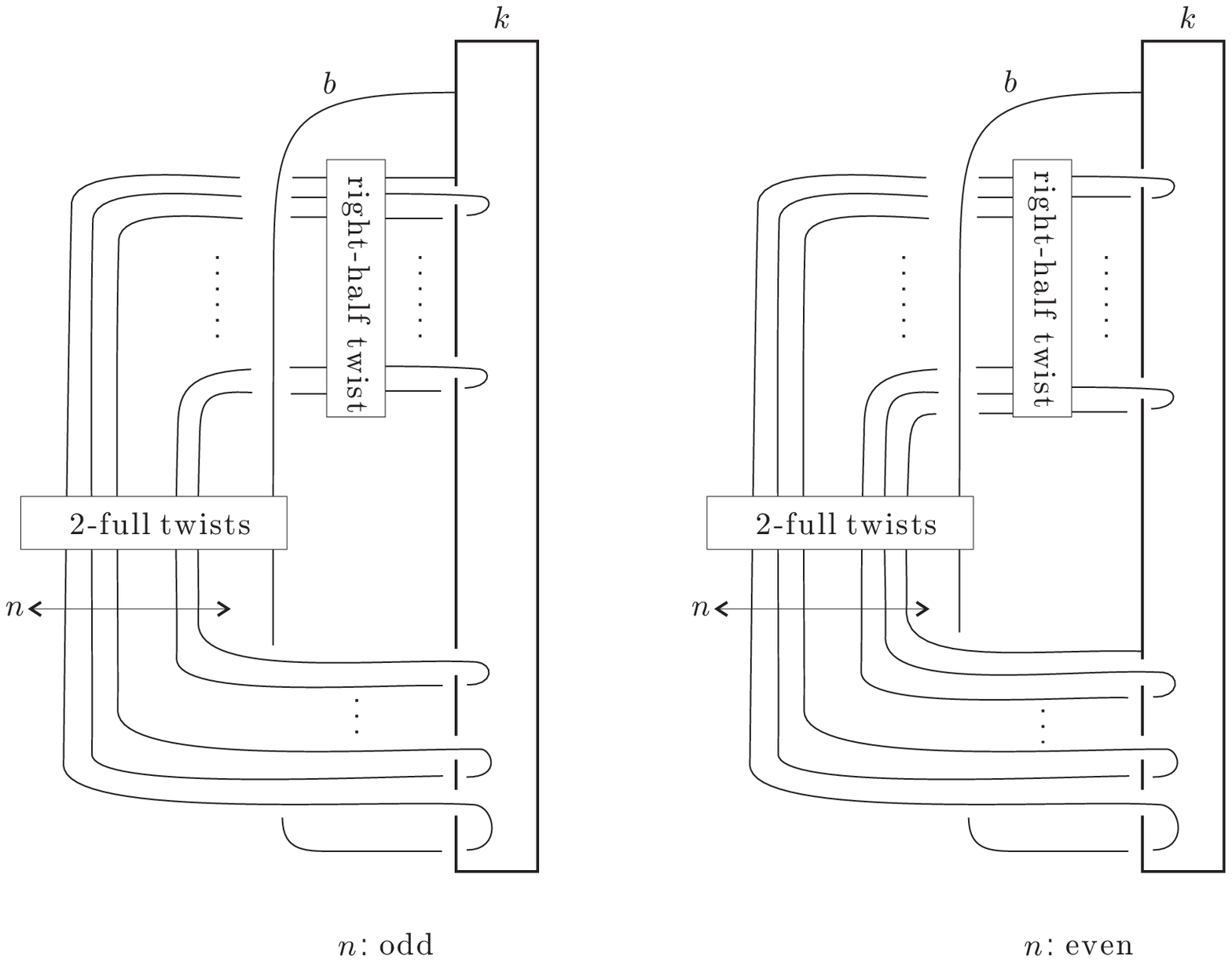}}
\vspace*{8pt}
\caption{The unknot $k$ with band $b$ (cont'd)}\label{fig:deform}
\end{figure}

\begin{figure}[th]
\centerline{
\includegraphics*[scale=0.6]{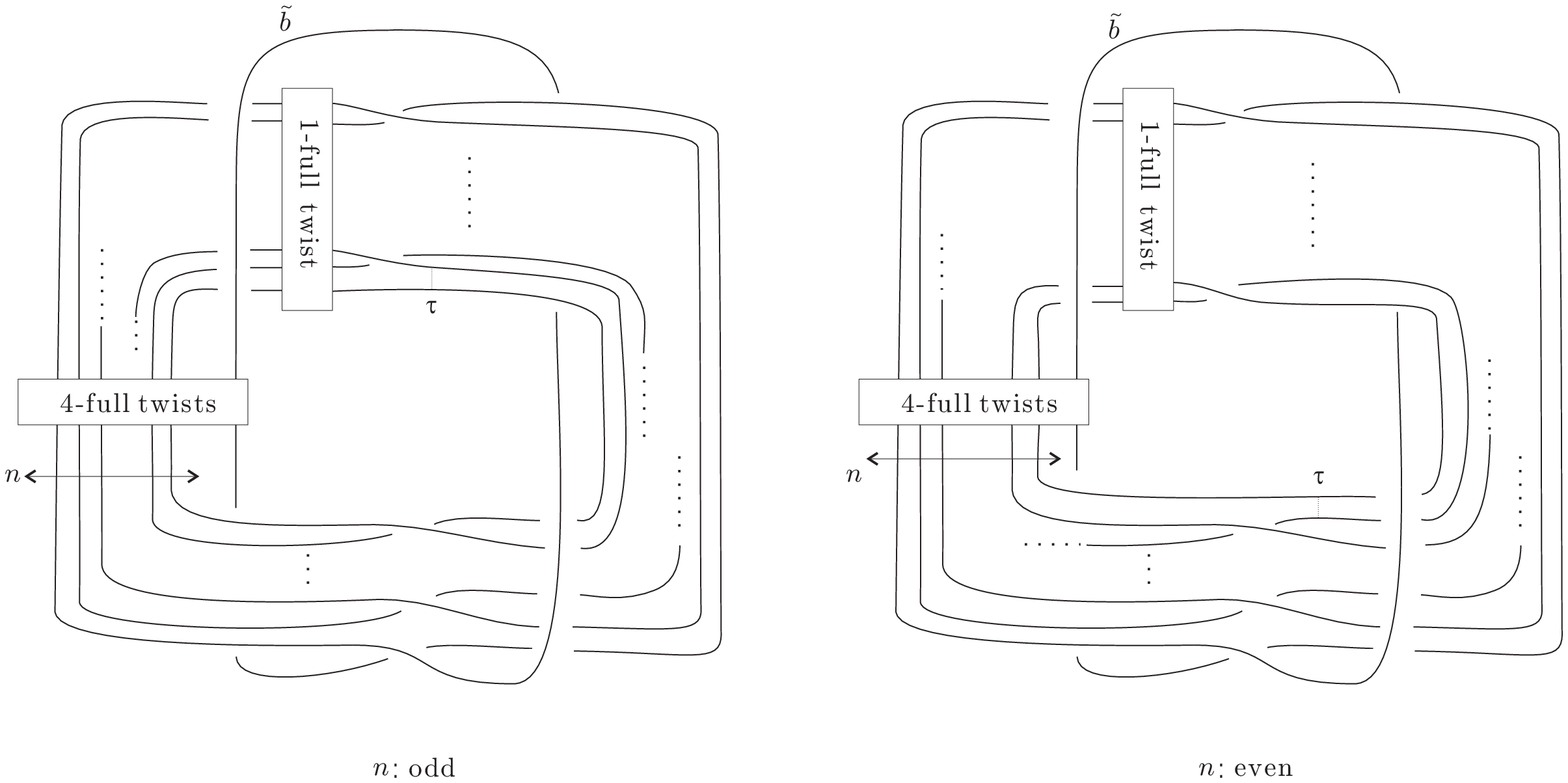}}
\vspace*{8pt}
\caption{$K_n$}\label{fig:final}
\end{figure}

First, we give an explicit description of $K_n$.
Consider the $1/0$-filling $T_n(1/0)$ of the tangle $T_n$.
Let $k$ denote the unknot $T_n(1/0)$.
To keep track of the framing, we indicate a band $b$ as shown in Fig.~\ref{fig:start}.
Since $k$ is trivial, it can be deformed so that it looks like a standard circle.
During this deformation, the band $b$ gets a complicated appearance.
In particular, we should be careful with the twists on the band.
See Fig.~\ref{fig:start} and Fig.~\ref{fig:deform}, where
a full twist means a right-handed full twist.
(In Fig.~\ref{fig:deform}, we indicate only the core of $b$ for simplicity.)
Let $\tilde{b}$ be the lift of $b$ in the double cover $S^3$ branched over $k$.
Then the core of $\tilde{b}$ is exactly the knot $K_n$, and
the framing of $\tilde{b}$ represents the slope corresponding to $0$-filling for $T_n$.
Figure \ref{fig:final} shows $K_n$, according to the parity of $n$.
In Fig.~\ref{fig:final}, $K_n$ has writhe $-3n$ and 
$\tilde{b}$ is represented as a flat band with $(-3)$-full twists before adding the $4$-full twists and the $1$-full twist indicated there as boxes.
Hence we see that $\tilde{b}$ has the framing $(-3n-3)+4(n+1)^2+n^2=5n^2+5n+1$ after performing those twists.
See \cite{EM,EM2} for this kind of procedure.

Hence we have the following.

\begin{proposition}\label{prop:slope}
The three toroidal slopes for $K_n$ correspond to $5n^2+5n+1,5n^2+5n+2,5n^2+5n+3$.
\end{proposition}

\begin{proposition}\label{prop:genus}
$K_n$ is a fibered knot of genus $(5n^2-n)/2$.
\end{proposition}

\begin{proof}
As seen in Fig.~\ref{fig:final}, $K_n$ is represented as a closed braid.
It is obvious that $K_n$ will be a closed positive braid after canceling negative crossings
with positive crossings coming from the $4$-full twists.
By \cite{St}, $K_n$ is fibered.
Moreover, its genus can be easily calculated by counting the number of crossings in the closed positive braid presentation,
because the Seifert surface obtained by applying Seifert's algorithm to the presentation is minimal genus \cite{St}.
\end{proof}

\noindent
\textbf{Proof of Theorem \ref{thm:main}.}
\enspace This immediately follows from Lemmas \ref{lem:3tor}, \ref{lem:tunnel}, \ref{lem:hyp},
and the fact that $K_n$'s are mutually distinct, which is a consequence of Proposition \ref{prop:genus}. \hfill $\Box$

\bigskip

\begin{figure}[th]
\centerline{
\includegraphics*[scale=0.5]{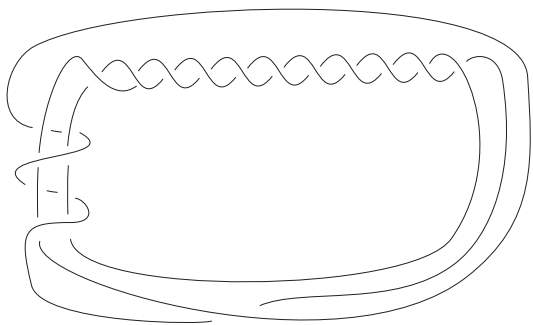}}
\vspace*{8pt}
\caption{$K_2$ with toroidal slopes $31,32,33$}\label{fig:k2}
\end{figure}

For reader's convenience, we exhibit the simplest $K_2$ of genus nine in Fig.~\ref{fig:k2}.

\section*{Acknowledgments}
The author would like to thank the referee for helpful suggestions.
Part of this work is supported by Japan Society for the Promotion of Science,
Grant-in-Aid for Scientific Research (C), 19540089.


\end{document}